\documentclass[12pt]{article}
\usepackage[cp866]{inputenc}
\usepackage[russian]{babel}
\usepackage{amssymb,amsmath,amsfonts}
\usepackage{srcltx}
\usepackage{graphicx}

\date{}
\newtheorem {theorem}{\bf Теорема}
\newtheorem {corol}{\bf Следствие}
\newtheorem {lemma}{\bf Лемма}

\newtheorem {stat}{\bf Утверждение}

\newcommand{\point}{\hspace{-1.75mm}{\bf.\ }}

\newcommand{\btheorem}{\begin{theorem}}
\newcommand{\etheorem}{\end{theorem}}

\newcommand{\bexam}{\begin{example}}
\newcommand{\eexam}{\end{example}}

\newcommand{\blemma}{\begin{lemma}\point}
\newcommand{\elemma}{\end{lemma}}

\newcommand{\bpro}{\begin{pro}\point}
\newcommand{\epro}{\end{pro}}

\newcommand{\bstat}{\begin{stat}\point}
\newcommand{\estat}{\end{stat}}

\newcommand{\bcorol}{\begin{corol}\point}
\newcommand{\ecorol}{\end{corol}}

\newcommand{\bnote}{\begin{note}\point}
\newcommand{\enote}{\end{note}}

\newcommand{\bconsec}{\begin{consec}\point}
\newcommand{\econsec}{\end{consec}}

\newcommand{\bdefin}{\begin{defin}\point}
\newcommand{\edefin}{\end{defin}}

\newtheorem {pro}{\bf Предложение}

\newcommand {\bproof }{{\par\medskip\noindent \bf Доказательство. }}
\newcommand {\eproof }{\hfill $\blacktriangle$ \\ \medskip}

\def\p0{\parindent 0pt}
\title{{\bf
\normalsize  И.~Ю.~Могильных\\
\large О продолжении пропелинейных структур кода
Нордстрома-Робинсона на код Хэмминга}
\thanks{Исследование
выполнено за счет гранта Российского научного фонда (проект
№14-11-00555)}}

\begin{document}

%\rm
%\pp
\maketitle

\begin{abstract}
Пропелинейным называется код, чья группа автоморфизмов содержит
подгруппу, регулярно действующую на коде. Такую подгруппу называют
пропелинейной структурой на коде. С помощью ПК перечислены
пропелинейные структуры на коде Норстрома-Робинсона и рассмотрена
задача их продолжения до пропелинейных структур на расширенном
коде Хэмминга длины 16. Последний результат основан на описании
разбиений кода Хэмминга длины 16 на коды Нордстрома-Робинсона
через плоскости Фано, которое приводится в данной работе. В итоге
получено рекордное число пропелинейных структур в классе
расширенных совершенных кодов длины 16.
\end{abstract}

\section{Введение}

Рассмотрим векторное пространство $F_2^n$ размерности $n$ над
полем из двух элементов. {\it Расстояние Хэмминга} между двумя
векторами из $F_2^n$ определяется как число координат, в которых
они различаются. Множество $C,$ $C\subset F_2^n$, называется {\it
двоичным кодом} с параметрами $(n,M,d)$,
 где $n$ - длина кода,
 $|C|=M$ и минимальное
 расстояние между различными кодовыми словами из $C$ равно $d$.
Код длины $n$, содержащий нулевой вектор $0^{n}$, будем называть
{\it приведенным}.

 Под группой автоморфизмов $Aut(F_2^n)$ пространства $F_2^n$ с метрикой Хэмминга будем понимать
 группу его изометрий относительно операции композиции. Хорошо известно, что $Aut(F_2^n)$
 исчерпывается преобразованиями вида $(x,\pi)$, где $\pi\in S_n$,
 $x\in F_2^n$, действующими следующим образом:

 $$(x,\pi)(y)=x+\pi(y)=(x_1+y_{\pi(1)},\ldots,x_n+y_{\pi(n)}).$$

В группе $Aut(F_2^n)$ выделим подгруппу $Sym(F_2^n)$, именуемую
{\it группой симметрий}:

$$Sym(F_2^n)=\{(0^n,\pi): \pi \in S_n\}.$$

{\it Группой автоморфизмов} (симметрий) кода $C$, обозначаемой
$Aut(C)$ ($Sym(C)$) называется стабилизатор множества его кодовых
слов в группе автоморфизмов (симметрий) пространства $F_2^n$.

Код Норстрома-Робинсона имеет 256 кодовых слов, длину $16$ и
является оптимальным с кодовым расстоянием 6 \cite{NR}. Известно,
что код с такими параметрами единственен с точностью до
автоморфизма $Aut(F_2^n)$. В этой связи далее кодом
Нордстрома-Робинсона будем называть любой код с его параметрами.
 Известно, что код Нордстрома-Робинсона является
единственным кодом, принадлежащим одновременно серии оригинальных
кодов Препараты и оригинальных кодов Кердока.

Среди кодов отдельно можно выделить {\it транзитивные} коды -
группа автоморфизмов которых содержит  подгруппу, действующую
транзитивно на кодовых словах. Если, более того, подгруппа
является {\it регулярной} (то есть ее порядок совпадает с
мощностью кода), то такой код называется {\it пропелинейным кодом}
\cite{Rifa1}, а сама подгруппа называется {\it пропелинейной
структурой} \cite{BMRS}.  Частным случаем пропелинейных кодов
являются $Z_4$-линейные коды. $Z_4$-линейность кода
Нордстрома-Робинсона была показана в работе \cite{FST}.

Напомним, что {\it кодом Препараты} называется максимальный
двоичный код длины $n=2^{m}$ для четного $m\geq4$ с кодовым
расстоянием 6. Первые представители этой серии кодов были впервые
построены в работе \cite{Prep}. Этот класс кодов обладает целым
рядом интересных свойств, связывающих его с классом совершенных
кодов. В частности, известно, что всякий код Препараты содержится
в единственном расширенном совершенном коде \cite{ZZS71}. Более
того, всякий оригинальный код Препараты \cite{Prep}, и в
частности, код Нордстрома-Робинсона, содержится в расширенном коде
Хэмминга.

\begin{stat}\label{u0}\cite{ZZS71}
Линейная оболочка кода Нордстрома-Робинсона совпадает с
расширенным кодом Хэмминга $H_{16}$ длины 16.
\end{stat}

Известно, что кодовые слова минимального ненулевого веса в коде
Хэмминга образуют систему четверок Штейнера. Как и в случае
расширенных совершенных кодов, кодовые слова минимального
ненулевого веса в приведенном коде Препараты образуют $3-$схему, а
кодовые слова веса 4 в приведенном коде Препараты, сдвинутом на
вектор веса 4 образуют $2-(n,4,1)$-схему \cite{SZZ71}.
Следовательно, всякое разбиение расширенного совершенного кода на
коды Препараты порождает разбиение системы четверок Штейнера этого
кода на $2-(n,4,1)$-схемы. Данное свойство для системы четверок
Штейнера называется {\it 2-разрешимостью}.

Отметим, что все известные на данный момент коды Препараты
\cite{Prep}, \cite{HKCS}, \cite{D}, \cite{BLW} являются
пропелинейными. Пропелинейность $Z_4$-линейных кодов следует из
конструкции \cite{HKCS}, более детальное исследование этих кодов и
их свойств можно найти в работе \cite{BPRZ}. Пропелинейность
неаддитивных кодов Препараты (конструкция Думера-Бакера) была
показана в \cite{RP}, также для этих кодов недавно показано
существование новой пропелинейной структуры \cite{ZinZin}.

Известно, что код Хэмминга имеет самые большие группы
автоморфизмов и симметрий в классе совершенных кодов \cite{ST1}.
Поэтому естественно предположить, что этот код имеет максимальное
число пропелинейных структур в классе совершенных кодов. Однако
вычислительных возможностей оказывается недостаточно для
верификации гипотезы даже для длины $n=16$. В этой связи
предлагается продолжать пропелинейные структуры с маломощных
пропелинейных подкодов (код Нордстрома-Робинсона, код Рида-Маллера
$RM(1,4)$) до пропелинейных структур на коде Хэмминга. В разделе 4
данной работы с помощью ПК показано, что для длины $n=16$ код
Хэмминга обладает самым большим количеством различных
пропелинейных структур среди пропелинейных совершенных кодов с
описанными пропелинейными структурами \cite{BMRS},\cite{MS}.
Результат существенно опирается на описание разбиений расширенного
кода Хэмминга длины 16 на коды Нордстрома-Робинсона, которое
приводится в разделе 2. О классификации пропелинейных структур на
коде Нордстрома-Робинсона речь идет в разделе 3.

\section{Разбиение расширенного кода Хэмминга на коды Нордстрома-Робинсона}

Вначале напомним некоторые известные факты, касающиеся групп
симметрий кодов Нордстрома-Робинсона и кодов Хэмминга.

 Через
$GL_m(F_q)$ будем обозначать группу невырожденных матриц порядка
$m$ над полем $F_q$ относительно умножения. Знакопеременная группа
$A_n$ состоит из четных подстановок на $n$ элементах. Известно,
что при $n\geq 5$ группа $A_n$ является {\it простой}, то есть не
содержит нетривиальных нормальных подгрупп.

Известно, что расширенный код Хэмминга $H_{2^m}$ является кодом
Рида-Маллера $RM(m-2,m)$, а код дуальный к нему - кодом
Рида-Маллера $RM(1,m)$. Напомним, что {\it кодом Рида-Маллера}
$RM(r,m)$ называется код, состоящий из векторов значений булевых
функций от $m$ переменных, представимых полиномами Жегалкина
степени не более $r$. Группа симметрий этих кодов исчерпывается
аффинными преобразованиями \cite{MS}. Если $f(v)$ - булевая
функция степени не более $r$, то для всякой невырожденной матрицы
$A$ порядка $m$ и вектора $b\in F_2^m$ функция $f(Av+b)$ имеет ту
же степень что и $f$, а подстановка координат, индуцируемая
заменой переменных $v\rightarrow Av+b$ принадлежит группе
симметрий кода. Отметим, что булева функция $f(v+b)$ отличается от
$f(v)$ лишь аффинной частью, другими словами, замена переменных
$v\rightarrow v+b$ оставляет неподвижным всякий класс смежности по
коду $RM(1,m)$.

\begin{stat} \label{u1} \cite{MS}
$$Sym(H_{2^m})= G\rightthreetimes H,$$ где $G \simeq  GL_m(F_2)$, $H \simeq Z_2^m$, при этом $\forall h \in
H$, $\forall g \in G$ выполнено

$$ h(x+RM(1,m))=x+RM(1,m), $$

$$g(x+RM(1,m))=y+RM(1,m) \mbox{ для некоторого } y\in F_2^m .$$
\end{stat}
Учитывая, что знакопеременная группа $A_8$ изоморфна $GL_4(F_2)$,
имеем следующее.
\begin{corol}\label{sl01} \cite{B07}
Группа симметрий расширенного кода Хэмминга $H_{16}$ изоморфна
$A_8\rightthreetimes Z_2^4$.
\end{corol}

 Итак, в группе $GL_4(F_2)\simeq A_8$ можно выделить подгруппу
$A_7$, порядок которой в 8 раз меньше порядка $A_8$. Эта подгруппа
отвечает нетривиальным симметриям кода Нордстрома-Робинсона,
мощность которого в 8 раз меньше мощности кода Хэмминга.

\begin{stat} \cite{B07} \label{u2} Пусть $N$-приведенный код
Нордстрома-Робин\-сона, содержащийся в коде $H_{16}$. Тогда
$Sym(N)\simeq A_7 \rightthreetimes Z_2^4< Sym(H_{16})$.

\end{stat}
Хорошо известно, что ядро всякого оригинального кода Кердока есть
код Рида-Маллера первого порядка \cite{MS}, поэтому ядро кода
Нордстрома-Робинсона есть код Рида-Маллера $RM(1,4)$ размерности
5.
 Рассмотрим разбиение кода Нордстрома-Робинсона на классы
смежности по его ядру: $$N=Ker(N)\cup \bigcup_{i=1\ldots7}
a_i+Ker(N).$$ Несложно видеть (см., например \cite{BMRS}), что
всякая симметрия кода оставляет на месте его ядро. Более того, в
силу Утверждения \ref{u1} подгруппа $Sym(N)$, изоморфная $Z_2^4$,
оставляет на месте каждый класс смежности $a+Ker(N)$, $a \in
H_{16}$, а группа, изоморфная $A_7$, действует на классах
смежности $a_i+Ker(N)$, $i=1,\ldots,7$. Отметим, что в силу того,
что $A_7$ - простая группа, никакие два ее элемента не будут
одинаково действовать на конечном множестве, то есть $A_7$ будет
действовать в привычном (естественном) смысле на 7 классах
смежности. Более того, действие является 2-транзитивным, так как
$A_7$ содержит подгруппу $GL_3(F_2)$.

\begin{corol} \label{sl2}
Пусть $N$-приведенный код Нордстрома-Робинсона. Тогда найдется
подгруппа группы $Sym(N)$, изоморфная $A_7$, которая действует
2-транзитивно на классах
 смежности $a_i+Ker(N), i=1,\ldots,7$ и оставляет на месте $Ker(N)$.
\end{corol}

\begin{lemma}\label{L1} Код Нордстрома-Робинсона, содержащийся в расширенном коде
Хэмминга $H_{16}$, является единственным с точностью до сдвига на
вектор из $H_{16}$. В коде $H_{16}$ существует ровно 8 приведенных
кодов Нордстрома-Робинсона, причем они единственны с точностью до
подстановки из $Sym(H_{16})$.
\end{lemma}

\bproof Пусть $N$, $N'$ - два кода Нордстрома-Робинсона,
содержащиеся в коде $H_{16}$ с кодовым словом $0^n$. С учетом
того, что коды Нордстрома-Робинсона единственны с точностью до
автоморфизма, $N$ и $N'$ получаются друг из друга подстановкой
$\pi\in Sym(H)$: $N=\pi(N')$. Следовательно, группы $Sym(N)$ и
$Sym(N')$ сопряжены при помощи подстановки $\pi\in S_{16}$:
$Sym(N)=\pi Sym(N') \pi^{-1}$.

 Таким
образом, число приведенных кодов Нордстрома-Робинсона в коде
Хэмминга равно числу подгрупп группы $Sym(H)$, сопряженных с
$Sym(N)$. С учетом Утверждения \ref{u2} и Следствия \ref{sl01},
число таких подгрупп равно числу подгрупп $A_8$, сопряженных с
$A_7$. В свою очередь, число подгрупп $A_8$, сопряженных с $A_7$,
равно $|A_8|/|\{g\in A_8: gA_7g^{-1}=A_7\}|$. Так как $A_7$ --
максимальна в $A_8$, то группа $\{g\in A_8: gA_7=A_7g\}$ либо
совпадает с $A_7$ либо с $A_8$. Последний случай невозможен, так
как знакопеременная группа не содержит нормальных подгрупп.

Итак, приведенных кодов Нордстрома-Робинсона в коде Хэмминга ровно
8. Заметим, что один приведенный код $N$ дает еще 7 различных
приведенных кодов: это коды $a_i+N$, где $a_i$, $i=1,\ldots,7$
есть представители классов смежности кода $N$ мощности $2^8$ по
его ядру мощности $2^5$.

 В силу транзитивности кода $N$ для любого $i$ найдется подстановка
$\pi_i$: $a_i+N=\pi_i(N)$. Другими словами, приведенные коды
Нордстрома-Робинсона, содержащиеся в коде Хэмминга, единственны с
точностью до подстановки. В силу того что линейная оболочка кодов
$N$ и $\pi_i(N)$ есть $H_{16}$, мы имеем $\pi_i\in Sym(H_{16})$.

 \eproof

\begin{lemma}\label{l2} Пусть $N$--код Нордстрома-Робинсона, $N=Ker(N)\cup \bigcup_{i=1\ldots7}
a_i+Ker(N).$ Тогда $\sum_{i=1,\ldots,7}a_i\in Ker(N)$.
\end{lemma}
\bproof Предположим что $\sum_{i\in \{1,\ldots,7\}}a_i\notin
Ker(N)$. Покажем что код, натянутый на вектора кода
Нордстрома-Робинсона, не является кодом Хэмминга.

%Предположим, что $$\sum_{i\in \{1,\ldots,7\}}a_i\notin \sum_{i\in
%L}a_i+Ker(N)$$ для любого непустого подмножества $L$ множества
%$\{1,\ldots,7\}$. Это равносильно тому,
% что для любого
%подмножества $I$ множества $\{1,\ldots,7\}$ имеет место
%$\sum_{i\in I}a_i\notin Ker(N)$, что в свою очередь, равносильно
%тому,
% что для различных подмножеств $I$ множества $\{1,\ldots,7\}$
%векторы $\sum_{i\in I}a_i$ принадлежат различным классам смежности
%по $Ker(C)$, откуда получаем, что

 Заметим, что $\sum_{i\in \{1,\ldots,7\}}a_i\in \sum_{i\in
L}a_i+Ker(N)$ для некоторого собственного подмножества $L$
множества $\{1,\ldots,7\}$. В противном случае получаем, что для
любых подмножеств $I, J$ множества $\{1,\ldots,7\}$ классы
$\sum_{i\in I}a_i+Ker(N)$ и $\sum_{i\in J}a_i+Ker(N)$ различны, и
 пространство, натянутое на кодовые слова кода $N$, имеет
размерность $12$, что противоречит Утверждению \ref{u0}.

Итак, $$\sum_{i\in \{1,\ldots,7\}}a_i\in \sum_{i\in L}a_i+Ker(N)$$
для некоторого непустого собственного подмножества $L$ множества
$\{1,\ldots,7\}$. Это равносильно тому, что $$a_i\in \sum_{k\in
I}a_k+Ker(N)$$ для некоторого подмножества $I$ множества
$\{1,\ldots,7\}$ мощности не более 5 и $i\in
\{1,\ldots,7\}\setminus I$. Рассмотрим произвольное $j\in
\{1,\ldots,7\}\setminus (i\cup I)$. Тогда для $a_j$ выполняется
аналогичное включение. Действительно, в силу 2-транзитивности
$Sym(N)$ (см. Следствие \ref{sl2}) на классах смежности
$a_i+Ker(N)$, $i=1,\ldots,7$ найдется подмножество $J$, $i, j
\notin J$: $a_j\in \sum_{k\in J}a_k+Ker(N)$. С одной стороны,
пространство, натянутое на векторы из $a_k+Ker(N),k \in
\{1,\ldots,7\}\setminus \{i,j\}$ имеет размерность не более 10, а
с другой стороны оно совпадает с пространством, натянутым на все
векторы кода Нордстрома-Робинсона, так как векторы $a_i$ и $a_j$
линейно выражаются через эти векторы. Получаем противоречие
Утверждению \ref{u0}. \eproof

Принимая во внимание предыдущую лемму и Утверждение \ref{u0},
получим следующее представление для кода Хэмминга:
\begin{theorem}\label{Hamrep}
Пусть $N$-- приведенный код Нордстрома-Робинсона, содержащийся в
коде Хэмминга $H_{16}$, $N=Ker(N)\cup \bigcup_{i=1\ldots7}
a_i+Ker(N)$. Тогда

$$H_{16}=Ker(N) \cup \bigcup_{i\in\{1,\ldots,7\}} (a_i+Ker(N)) \cup $$ $$\bigcup_{\{i,j\}\subset \{1,\ldots,7\}} (a_i+a_j+Ker(N))\cup \bigcup_{\{i,j,k\}\subset\{1,\ldots,7\}}(a_i+a_j+a_k+Ker(N)).$$
\end{theorem}

Теперь рассмотрим случай, когда пара кодов Нордстрома-Робинсона,
содержащихся в одном коде Хэмминга, не пересекаются. Согласно
Лемме \ref{l2}, один код из другого получаются сдвигом на вектор
из $H_{16}$.

\begin{lemma}\label{linter}
Коды  Нордстрома-Робинсона $a+N$ и $b+N$, содержащиеся в коде
Хэмминга $H_{16}$ не пересекаются тогда и только тогда, когда
$a+b\in a_i+a_j+a_k+Ker(N)$ для некоторого трехэлементного
подмножества $\{i,j,k\}$ множества $\{1,\ldots,7\}$.
\end{lemma}
\bproof Утверждение достаточно доказать в случае $a=0$, а вектор
$b$, согласно Теореме \ref{Hamrep}, лежит в $a_i+Ker(N)$,
$a_i+a_j+Ker(N)$ или $a_i+a_j+a_k+Ker(N)$. Однако в первых двух
случаях коды $N$ и $b+N$ содержат одновременно класс $a_i+Ker(N)$.
С другой стороны, в  третьем случае $b+N$ состоит из классов
смежности по $Ker(N)$ с представителями $a_i+a_j+a_k+a_l$, где
$l\in \{1,\ldots,7\}$. В трех случаях эти представители являются
суммой двух различных представителей классов смежности $N$ по
$Ker(N)$, а в четырех случаях, принимая во внимание Лемму \ref{l2}
- суммой трех различных представителей классов смежности $N$ по
$Ker(N)$. Следовательно, по Теореме \ref{Hamrep} код $N$ не
пересекается с кодом $N+b$. \eproof

\begin{theorem}\label{tpart}
Пусть $N$-- приведенный код Нордстрома-Робинсона, $N=Ker(N)\cup
\bigcup_{i=1\ldots7} a_i+Ker(N),$ $H_{16}$ - код Хэмминга,
содержащий код $N$. Для всякой плоскости Фано $S$ на точках
 \{1,\ldots,7\}
$$N\cup\bigcup_{\{i,j,k\}\in S} (a_i+a_j+a_k+N)$$
%$$N\cup\{a_i+a_j+a_k+N:\{i,j,k\}\in S \}$$

есть разбиение кода Хэмминга $H_{16}$ на коды Нордстрома-Робинсона
и других разбиений, содержащих $N$, нет.

\end{theorem}

\bproof Рассмотрим произвольное разбиение $$N\cup \bigcup_{b\in
B}b+N$$ на коды Нордстрома-Робинсона. Согласно Лемме \ref{linter},
для всякого $b\in B$ имеет место: $$ b=a_{i_b}+a_{j_b}+a_{k_b},$$
для некоторого трехэлементного подмножества $\{i_b,j_b,k_b\}$
множества $\{1,\ldots,7\}$.

 Рассмотрим пересечение
множеств $\{i_b,j_b,k_b\}$ и $\{i_{b'},j_{b'},k_{b'}\}$ для $b,
b'\in B$. Если эти множества пересекаются по 2 элемента, то $b+b'$
есть сумма двух представителей классов смежности кода $N$ по
$Ker(N)$ и согласно Лемме \ref{linter} коды $b+N$ и $b'+N$
пересекаются и не могут входить в одно разбиение. Если множества
не пересекаются, то согласно Лемме \ref{l2}, кодовое слово $b+b'$
принадлежит коду $N$, что влечет что $b+N$ и $b'+N$ имеют общие
кодовые слова. Наконец, в случае, когда множества
$\{i_b,j_b,k_b\}$ и $\{i_{b'},j_{b'},k_{b'}\}$ имеют один общий
элемент, вектор $b+b'$, согласно Лемме \ref{l2}, есть сумма
представителей трех различных классов смежности кода $N$ по
$Ker(N)$ и коды $b+N$ и $b'+N$ не пересекаются. Утверждение
теоремы вытекает из того что любая пара из 7 прямых плоскости Фано
пересекается в одной точке. \eproof

Два разбиения $\bigcup_{i\in I}C_i$ и $\bigcup_{j\in J}D_j$ кода
Хэмминга $H$ назовем {\it изоморфными}, если существует $\pi \in
Sym(H)$, переводящая элементы одного разбиения в другое. Заметим,
что $\pi\in Sym(C_i)$, где $C_i$ приведенный элемент разбиения.

\begin{corol}
Существует ровно два неизоморфных разбиения кода Хэмминга на коды
Нордстрома-Робинсона.
\end{corol}

\bproof Рассмотрим пару произвольных разбиений кода Хэмминга на
коды Нордстрома-Робинсона. В силу Леммы \ref{L1} без ограничения
общности их можно полагать содержащими один и тот же приведенный
код Нордстрома-Робинсона $N$.

 Согласно Теореме \ref{tpart} два разбиения на коды
 Нордстрома-Робинсона различаются выбором плоскостей Фано и имеют следующий вид:
$$N\cup\bigcup_{\{i,j,k\}\in S} (a_i+a_j+a_k+N),$$
$$N\cup\bigcup_{\{i,j,k\}\in S'} (a_i+a_j+a_k+N),$$
где $S$, $S'$ -- множество прямых двух плоскостей Фано на точках
$\{1,\ldots,7\}$. Принимая во внимание Утверждение \ref{u1},
изоморфизм двух разбиений может быть осуществлен лишь при помощи
подстановки из $A_7$, существование которой равносильно
изоморфизму плоскостей Фано, соответствующих этим разбиениям.
 В заключение заметим, что группа $A_7$ имеет две
орбиты на всех плоскостях Фано \cite{B07} на фиксированном
7-элементном множестве. При этом каждая орбита состоит из 15
плоскостей, а плоскости одной орбиты получаются из другой с
помощью свитчинга Паш-конфигурации. \eproof

\section{Пропелинейные структуры на коде Нордстрома-Робинсона}
Напомним, что {\it графом минимальных расстояний}
 кода называется граф, вершинами которого являются
кодовые слова, соединенные ребрами, если расстояние Хэмминга между
ними есть кодовое расстояние. Известно, что группа автоморфизмов
всякого кода Препараты (в том числе кода Нордстрома-Робинсона)
изоморфна группе автоморфизмов его графа минимальных расстояний
(при $n\geq 2^{12}$ следует из \cite{SA03}, \cite{Mog09}, для
произвольного $n$ из работы \cite{PhFer}). Аналогичное свойство
выполнено для совершенных кодов и расширенных совершенных кодов
\cite{AVG98}, \cite{OST08} при $n\geq 16$.
 По этой причине поиск
пропелинейных структур на вышеперечисленных классах кодов или
проверку транзитивности  наиболее вычислительно быстро можно
реализовать через перечисление регулярных или поиск транзитивных
подгрупп группы автоморфизмов графа минимальных расстояний этих
кодов с использованием пакета MAGMA. В теореме ниже под
сопряженными и изоморфными пропелинейными структурами понимаются
аналогичные понятия для групп.
 Заметим, что для всякой пропелинейной структуры $(C,*)$ на коде
 $C$ выполнено (см. \cite{BMRS}):

 $$ |\{\pi:(x,\pi_x)\in (C,*)\}|\leq |C|/|Ker(C)|, $$
 а коды, достигающие этой границы, называются
{\it нормализованно-пропелинейными}.

\begin{theorem} На коде Нордстрома-Робинсона существует ровно 338
классов сопряженности пропелинейных структур, разбивающихся на 250
классов изоморфизма. Существует 28 нормализованно-пропелинейных
классов сопряженности, которые разбиваются на 25 классов
изоморфизма.
\end{theorem}

%{\it Замечание.} Отметим, что на коде Нордстрома-Робинсона
%существуют нормализованно-пропелинейные структуры, изоморфные
%ненормализованно-пропелинейным. В этой связи естественно ставить
%вопрос о подходящих инвариантах для классов изоморфизма
%пропелинейных структур на одном коде.

\section{Продолжение пропелинейных структур
 на коде \- Нордстрома-Робинсона на код Хэмминга}

Для пропелинейной структуры $(C,*)$ введем следующее обозначение
$\Pi_C=\{\pi:(x,\pi) \in (C,*)\}$.
 Пусть $C$ является подкодом кода $D$. Пропелинейную структуру
$(D,*)$ на коде $D$ будем называть {\it продолжением}
пропелинейной структуры $(C,*)$ на коде $D$, если $(C,*)<(D,*)$.
Продолжение будем называть {\it  в узком смысле}, если
$\Pi_C=\Pi_D$. В данной работе рассматриваются лишь продолжения
пропелинейных структур в узком смысле.

Частным случаем таких продолжений будут $Z_4$-линейные совершенные
коды. Рассмотрим произвольную подгруппу $G$ группы $Z_4^{n/2}$.
Каждому элементу этой группы поставим в соответствие двоичный
вектор $F_2^n$, произведя замены с помощью {\it отображения Грэя}
$\phi$: $0\rightarrow 00, 1\rightarrow 01, 2\rightarrow 11,
3\rightarrow 01.$ Образ подгруппы $G$ под действием отображения
Грэя называется $Z_4$-{\it линейным кодом}. Рассмотрим кодовое
слово $x$ некоторого $Z_4$-линейного кода $D$. Зададим
$Z_4$-линейную структуру на коде $C$. Для этого назначим каждому
$x$ подстановку $\pi_x$ по правилу: если
$(\phi^{-1}(x))_i=0(mod\mbox{ }2)$, то $\pi_x(2i)=2i,
\pi_x(2i+1)=2i+1$, иначе $\pi_x(2i)=2i+1, \pi_x(2i+1)=2i$.
Очевидно, что группа $\{(x,\pi_x):x \in D\}\simeq G$ будет
регулярно действовать на $D$.

\begin{stat}\label{prper} Пусть $(P,*)$ является $Z_4$-линейной
структурой на коде Препараты, тогда она продолжается в узком
смысле до $Z_4$-линейного расширенного совершенного кода $C$,
содержащего $P$.
\end{stat}
\bproof В силу сказанного выше, для любых двух кодовых слов $x,y$
некоторого $Z_4$-линейного кода имеет место: $\pi_x=\pi_y
\Leftrightarrow 2\phi^{-1}(x)=2\phi^{-1}(y)$. Следовательно,
$Z_4$-линейная структура $(P,*)$ продолжается в узком смысле до
$Z_4$-линейной структуре на надкоде $(C,*)$ тогда и только тогда
когда $2\phi^{-1}(P)=2\phi^{-1}(C)$, что выполняется для
$Z_4$-линейного кода Препараты и совершенного кода $C$ его
содержащего (см. \cite{BPRZ}, Следствие 3.11).

\eproof

\begin{stat}\label{exten} Если пропелинейная структура $(C,*)$ продолжается в
узком смысле до пропелинейной структуры $(D,*)$, тогда

1. $(D,*)\simeq (C,*)\rightthreetimes  Z_2^{log_2(|D|/|C|)}$

2. Существует разбиение кода $D$ на сдвиги кода $C$, которые
совпадают с левыми классами смежности группы $(D,*)$ по $(C,*)$.

\end{stat}
\bproof Заметим, что в силу условия $\Pi_C=\Pi_D$, в качестве
представителя любого левого смежного класса $x*C$ можно выбрать
$x$, которому в пропелинейной структуре $(D,*)$ назначен
автоморфизм с тождественной подстановкой: $(x,id)\in (D,*)$.
Рассмотрим следующие элементарные абелевы 2-группы $C_{id}=\{x\in
C: (x,id) \in (C,*) \}$ и ее надгруппу $D_{id}=\{x\in D: (x,id)
\in (D,*) \}$. Очевидно, что фактор-группа $D_{id}/C_{id}$
является нормальной подгруппой $(D,*)$, следовательно $(D,*)\simeq
(C,*)\rightthreetimes Z_2^{log_2(|D|/|C|)}$. Более того, в силу
замечания выше, каждый класс смежности по $(C,*)$ есть сдвиг кода
$C$.

\eproof

{\bf {\large Алгоритм 1.}}

{\bf Вход:} Пропелинейная структура $(N,*)$ на коде
Нордстрома-Робинсона $N$.

{\bf Шаг 1.} Перебор разбиения $N,c_1+N,\ldots,c_7+N$ кода
Хэмминга $H_{16}$ на коды Нордстрома-Робинсона с помощью
представления разбиений через плоскости Фано (Теорема 2).

{\bf Шаг 2.} Перебор порождающих для элементарной абелевой
2-группы $L$ порядка 8 среди классов смежности:
$c_1+Ker(N),\ldots,c_7+Ker(N)$.

{\bf Шаг 3.} Если группа, порожденная $L$ и $(N,*)$, имеет порядок
$|H_{16}|$, то согласно Предложению \ref{exten} пропелинейная
структура $(N,*)\rightthreetimes L$ на коде $H_{16}$ есть
продолжение структуры $(N,*)$ в узком смысле, иначе возвращение на
Шаг 2.

{\bf Выход:} Все пропелинейные структуры на коде Хэмминга
$H_{16}$, продолжающие $(N,*)$ в узком смысле.

С использованием Алгоритма 1 с помощью  ПК получена следующая
классификация. Приводимая в теореме 4 нижняя оценка числа
неизоморфных пропелинейных структур установлена с помощью порядка
централизаторов элементов пропелинейной структуры.

\begin{theorem}
Существует 3057 классов сопряженности и хотя бы 2284 классов
изоморфизма пропелинейных структур на расширенном совершенном коде
Хэмминга $H_{16}$, продолжающих пропелинейные структуры на коде
Нордстрома-Робинсона. Существует единственная пропелинейная
структура (с точность до сопряжения) на коде Нордстрома-Робинсона,
которая не имеет продолжения на $H_{16}$.
\end{theorem}

{\bf Заключение.} Проблема построения пропелинейных структур
(регулярных подгрупп группы автоморфизмов) в классе оптимальных
кодов является сложной проблемой. В данной работе предлагается
получать структуры, продолжая их с пропелинейных подкодов, чья
группа автоморфизмов естественно вкладывается в группу
автоморфизмов оптимального кода. Данный подход реализованный для
кодов Хэмминга и кодов Нордстрома-Робинсона позволил получить
рекордное количество пропелинейных структур в классе расширенных
совершенных кодов длины 16. Помимо естественного интереса с точки
зрения алгебро-комбинаторной теории кодирования, данный результат
может найти применение в криптографии в системах, сходных с
криптосистемой МакЭлиса \cite{Mac}.

%\begin{theorem}
%Пусть $N$--код Нордстрома-Робинсона, $N=Ker(N)\cup
%\bigcup_{i=1\ldots7} a_i+Ker(N),$ $H_{16}$ - код Хэмминга,
%содержащий код $N$. Тогда
% $$N\cup
%\bigcup_{\{i,j,k\}\in S}a_i+a_j+a_k+N$$ является разбиением
%$H_{16}$ на коды Нордстрома-Робинсона тогда и только тогда когда
%$S$-множество прямых некоторой плоскости Фано на множестве точек
%$\{1,\ldots,7\}$.
%\end{theorem}
%\bproof
%\eproof

Автор выражает благодарность Кротову Денису Станиславовичу и
Соловьевой Фаине Ивановне за полезные обсуждения в ходе работы над
статьей.

 \bigskip

Могильных Иван Юрьевич

Институт математики им. С.Л.Соболева СО РАН,
ivmog@math.nsc.ru
\end{document}